\documentclass[12pt,draftcls,onecolumn]{IEEEtran}
\usepackage{mathrsfs,amsmath,latexsym,amssymb}

\newtheorem{thm}{Theorem}[section]

\newtheorem{assum}{Assumption}[section]
\newtheorem{cor}{Corollary}[section]
\newcommand{\argmax}{\mathop{\mbox{\rm arg\,max}}}
\newcommand{\argmin}{\mathop{\mbox{\rm arg\,min}}}
\def\Bbb{\mathbb } \def\real {{\Bbb R}} 
\newcommand{\PI}{\mathop{\mbox{\small PI}}}
\newcommand{\app}{\mathop{\mbox{\footnotesize app}}}
\begin{document}
\sloppy

\title{Approximate Constrained Discounted Dynamic Programming with Uniform Feasibility and Optimality}
\author{Hyeong Soo Chang
\thanks{H.S. Chang is with the Department of Computer Science and Engineering at Sogang University, Seoul 121-742, Korea. (e-mail:hschang@sogang.ac.kr).}%
}

\maketitle
\begin{abstract}
We consider a dynamic programming (DP) approach to approximately solving
an infinite-horizon constrained Markov decision process (CMDP) problem with a fixed 
initial-state for the expected total discounted-reward criterion
with a uniform-feasibility constraint of the expected total discounted-cost
in a deterministic, history-independent, and stationary policy set.
We derive a DP-equation that recursively holds 
for a CMDP problem and its sub-CMDP problems, where
each problem, induced from the parameters of the original CMDP problem,
admits a uniformly-optimal feasible policy in its policy set associated 
with the inputs to the problem.
A policy constructed 
from the DP-equation is shown to achieve the optimal values, defined
for the CMDP problem the policy is a solution to, at all states.
Based on the result, we discuss off-line and on-line computational algorithms, 
motivated from policy iteration for MDPs, whose output sequences have 
local convergences for the original CMDP problem.
\end{abstract}


\section{Introduction}

Consider a problem of an infinite-horizon Markov decision process (MDP)
with the expected total discounted-reward criterion with a uniform-feasibility 
constraint of the expected total discounted-cost 
(see, e.g.,~\cite{altman}~\cite{feinberg}).
The constrained MDP (CMDP) problem is associated with a finite state-set $X$ and a finite action-set $A$ 
where $A(x)$ denotes the non-empty set of admissible actions in $A$ at $x$.
We define a policy $\pi$ as a mapping from $X$ to $A$ whose graph $G(\pi)$ is equal to $\{(x,\pi(x))| x\in X, a \in A(x)\}$ and $\Pi$ as 
the set of all possible policies. 
Thus any policy in $\Pi$ is \emph{deterministic (or pure), history-independent (or Markovian), and
stationary}.
When a fixed policy $\pi$ in $\Pi$ is followed in the CMDP, the underlying process 
acts as a Markov chain 
whose state-transition 
dynamics is as follows: once an action $\pi(x)$ is taken at $x$ in $X$, $x$ makes a random transition to $y$ in $X$ according to the probability specified by $P_{xy}^{\pi(x)}$. 
The system obtains
a reward of $R(x,\pi(x))$ and a cost of $C(x,\pi(x))$ where $R$ and $C$ are functions from 
$\{(x,a)|x\in X, a\in A(x)\}$ to $\real$, respectively.

Denote the set of all possible real-valued functions over a finite set $S$ as $B(S)$.
(In the sequel, any operator is applied componentwise for the elements in $B(S)$.)
We define the value functions $V^{\pi}$ and $J^{\pi}$ of the expected total discounted-reward and the expected total discounted-cost of $\pi$ in $\Pi$, respectively, in $B(X)$, such that the values of $\pi$ at (an initial
state) $x$ in $X$
as
\[
V^{\pi}(x) = E\left [\sum_{t=0}^{\infty} \gamma^t R(X_t,\pi(X_t)) \biggl |X_0=x \right ]
\]
and
\[
J^{\pi}(x) = E\left [\sum_{t=0}^{\infty} \beta^t C(X_t,\pi(X_t)) \biggl |X_0=x \right  ],
\]
respectively, where $X_t$ is a random variable representing a state at time $t$ by following $\pi$ and 
discounting factors $\gamma$ and $\beta$ are in $(0,1)$, respectively.
The expectation is taken over the probability distribution over all possible random trajectories of
an infinite length obtained by following $\pi$ with the starting state $x$.

Define a mapping $\Phi:\Pi \rightarrow 2^{\Pi}$ such that 
\[
\Phi(\pi) = \{g| g\in \Pi, J^g\leq J^{\pi} \}, \pi\in \Pi.
\]
We say that $g$ is \emph{uniformly $\pi$-feasible} if $g$ is in $\Phi(\pi)$.
Thus for each $\pi$, $\Phi(\pi)$ is the set of the policies that satisfy
\emph{uniform-feasibility} with respect to $\pi$ and is not empty because $\pi$ 
is in $\Phi(\pi)$.

Define the \emph{optimal value function} $V^*_c$ in $B(X\times \Pi)$ as
\[ V^*_c(x,\pi) = \max_{g\in \Phi(\pi)} V^g(x), x\in X, \pi\in \Pi.
\]

The problem we deal with in this paper, referred to as $M^c_{\pi^c}(x_0)=(X,A,P,R,C,\pi^c,x_0)$, is to find 
$\pi^*_c$ in $\Phi(\pi^c)$ such that $V^*_c(x_0,\pi^c) = V^{\pi^*_c}(x_0)$ for a given \emph{threshold} policy $\pi^c$ in $\Pi$ and a fixed initial-state $x_0$ in $X$.
Because $x_0$ is fixed and $\Phi(\pi^c)$ is non-empty, such $\pi^*_c$ always exists.

The CMDP problems are mostly set up in the literature
with a specific initial-state $x$ or with an initial-state probability distribution over a \emph{randomized} policy set (a policy's codomain corresponds to the set of the possible probability distributions over $A(x), x\in X$).
This is because
due to the constraint(s), the feasibility of a policy and the optimality might depend on the initial state~\cite{altman}. 
Seemingly impossible, no DP optimality-equation that covers certain optimal-substructure and overlapping subproblem characteristics has been established with this formulation.
Much attention has been paid to finding out relevant structures in best randomized policies (see, e.g.,~\cite{altman}~\cite{feinberg} and the references therein) to exploit those while devising solution algorithms.
In particular, because a CMDP with such formulation can be cast into an equivalent
linear programming (LP) problem (e.g., by ``occupation measure"~\cite{altman}), employing LP in
solution algorithms is a very common approach.

It should be noted that the problem setup also includes a threshold function in $B(X)$ as a parameter to identify the feasibility of a policy.
In $M^c_{\pi^c}(x_0)$, the threshold function is given by $J^{\pi^c}$.
We consider such threshold function determined by the value function of a threshold policy
in order to exploit a substructural relationship between CMDP problems in 
the DP-equation to be discussed below.
(If this is critical, one can choose a policy whose value function approximates the given threshold 
function. Refer to the further comment in the concluding remarks section.)

A closely related problem to $M^c_{\pi^c}(x_0)$ is to find
\[\argmax_{\pi\in \Phi_{\delta}(\pi^c)} \left ( \sum_{x\in X} \delta(x) V^{\pi}(x) \right ),
\] 
where $\delta$ is a fixed initial-state distribution and
\[\Phi_{\delta}(\pi^c) = \left \{g \biggl | g\in \Pi, \sum_{x\in X} \delta(x) J^g(x) \leq \sum_{x\in X} \delta(x) J^{\pi^c}(x) \right \}.
\] Feinberg showed that this is an NP-hard problem~\cite{feinberg} when the input size is determined by $|X|$, $|A|$, and the number of constraints.
The problem is no easier to solve than LP-solvable CMDPs.
Feinberg, in fact, proved that LP is not applicable. A non-linear and non-convex mathematical program (MP) is equivalent to this combinatorial problem in the sense that the MP is feasible if and only if $\Phi_{\delta}(\pi^c)$ is not empty (cf., P1 in~\cite[Theorem 3.1]{feinberg}).
Furthermore, a Lagrangian approach with a ``saddle point property"~\cite{altman}~\cite{kadota} would not work due to the \emph{deterministic} policy space as remarked in~\cite{changexact}.
No DP-equation is available yet to the problem either so that the title of~\cite{feinberg}, ``constrained discounted DP," is somewhat misleading.
A notable work for this problem is an exact iterative-algorithm that follows the idea of policy iteration (PI)~\cite{changexact} with a characterization of the entire feasible policy. 
In Section~\ref{sec:algo}, we discuss a PI-type algorithm that uses a similar idea for (approximately) solving $M^c_{\pi^c}(x_0)$.

Chen and Blankenship~\cite{chen} provided DP-equations that hold with value functions defined over $X\times B(X)$ for CMDPs with uniform-feasibility constraint(s) 
when an optimal policy is sought for over a deterministic but \emph{history-dependent} policy space with a fixed initial-state 
although no specific (iterative) computational algorithms were presented in that work.
Later,
Chen and Feinberg~\cite{chen2} proposed a value iteration (VI)-type algorithm, based on the results 
of the DP-equations in~\cite{chen}, over a deterministic but \emph{non-stationary} policy space.
The algorithm iteratively solves a sequence of finite-horizon MDPs with a sequence of increasing 
horizon sizes. 
The sequence of the solutions to the MDPs is shown to converge to an optimal policy as the horizon 
size approaches infinity.
Putting aside the soundness of the approach, computing their DP-equation is rather problematic
because the equation is involved with a maximization operator over $X\times B(X)$.
Another DP-approach was taken by Piunovskiy and Mao~\cite{piu} over a \emph{randomized} policy space.
A value function is defined as a function of a state distribution and the expected total cost
and a VI-type algorithm, based on the function, was presented where the complexity of computing the DP-equation would be again an issue.

The DP optimality-equation by Bellman for an MDP allows to construct
a \emph{uniformly-optimal} policy in $\Pi$
that 
achieves the optimal values, defined for the MDP, at \emph{all} states.
We first derive DP-equation that recursively holds for a CMDP problem and
its sub-CMDP problems, where each problem is induced from the parameters 
of $M^c_{\pi^c}(x_0)$ and admits a uniformly-optimal feasible policy in
its policy set associated with the problem. A policy constructed 
from the DP-equation is shown to achieve the optimal values, defined
for the CMDP problem the policy is a solution to, at all states.
In rough terms, $M^c_{\pi^c}(x_0)$ ``contains" an associative CMDP problem whose
solution policy is uniformly-feasible and uniformly-optimal and the value 
of the solution policy at $x$ is a lower bound to $V^*_c(x,\pi^c)$.
This provides an accompanying result to the theories of MDP that can
be used for exploiting a DP structure to approximately solve a CMDP problem
formulated with a fixed initial-state and the uniform-feasibility.
We then present computational algorithms, motivated from PI for MDPs,
within both off-line and on-line contexts whose output sequences have local-convergences for $M^c_{\pi^c}(x_0)$. The initialization of the algorithms
can use a solution to the associative CMDP of $M^c_{\pi^c}(x_0)$ solved by DP.

\section{Dynamic Programming Optimality Equation}
\label{sec:dp}

Define $T_{\phi}:B(X)\rightarrow B(X)$ and $U_{\phi}:B(X)\rightarrow B(X)$
for $\phi\in \Pi$ such that
$T_{\phi}(u)(x) = R(x,\phi(x)) +\gamma \sum_{y\in X} P_{xy}^{\phi(x)} u(y)$ 
and 
$U_{\phi}(u)(x) = C(x,\phi(x)) +\beta \sum_{y\in X} P_{xy}^{\phi(x)} u(y)$,
respectively, for any $u\in B(X)$ and $x\in X$. 

Suppose that we have a set $\{A^{\pi}, \pi \in \Pi \}$ whose element $A^{\pi}$ is a non-empty subset of $A$ such that $\pi(x) \in A^{\pi}(x)$ and $A^{\pi}(x) \subseteq A(x)$ for all $x$ in $X$.
Define $F:\Pi \rightarrow 2^{\Pi}$ such that 
\[
F(\pi) = \{g| g\in \Pi, g(x) \in A^{\pi}(x) \mbox{ } \forall x\in X\}, \pi\in \Pi.
\] It is the set of all policies whose codomains are equal to $A^{\pi}$.
If $F(\pi)$ is a subset of $\Phi(\pi)$, then $F(\pi)$ is called a \emph{$\pi$-feasible policy set induced from} $A^\pi$.
If for \emph{all} $\pi$ in $\Pi$, $A^\pi$ induces a $\pi$-feasible policy set, then 
$\{A^{\pi}, \pi \in \Pi \}$ is referred to as a \emph{DP-inducible set}.

We provide an example of the DP-inducible set $\{A^{\pi}, \pi \in \Pi \}$. Let
\begin{eqnarray} 
\label{eqn:DPset}
A^{\pi}(x) = \left \{a \biggl | a\in A(x), C(x,a) + \beta\sum_{y\in X}P_{xy}^a J^{\pi}(y) \leq J^{\pi}(x) \right \}, x\in X, \pi\in \Pi.
\end{eqnarray} Then because for any $\phi$ in $F(\pi)$, $U_{\phi}(J^{\pi}) \leq J^{\pi}$, we have that $J^{\phi}\leq J^{\pi}$ (see, e.g.,~\cite{ross}). 
It is important to remark that $A^{\pi}$'s across $\pi$ in $\Pi$ should have 
the same \emph{form} and the same \emph{domain} for each input parameter 
in order for any subproblem to be well-defined with a smaller input size 
than a problem that recursively calls the subproblem.
This is from the fundamental property of DP. A CMDP problem needs to exhibit
optimal substructure such that an optimal solution to the problem contains 
within it optimal solutions to sub-CMDP problems.
Throughout the paper, we will assume that the following assumption holds:

\begin{assum}
\label{ass:ext}
A DP-inducible set $\{A^{\pi}, \pi \in \Pi \}$ is given for $M^c_{\pi^c}(x_0)$. 
\end{assum}

The assumption allows us to define the \emph{uniformly-optimal} value function
$V^*$ in $B(X\times \Pi)$ as
\[ V^*(x,\pi) = \max_{g\in F(\pi)} V^g(x), x\in X, \pi\in \Pi.
\]
Therefore, by definition $V^*(x,\pi)$ is a lower bound to $V^*_c(x,\pi)$ for any $x$ in $X$ and $\pi$ in $\Pi$.

\begin{thm}
\label{thm:ext}
For any $\pi\in \Pi$, there exists a uniformly-optimal $\pi$-feasible
$\pi^*$ in $F(\pi)$ such that $V^*(x,\pi) = V^{\pi^*}(x)$ for all $x\in X$.
\end{thm}
\proof
This directly follows from Assumption~\ref{ass:ext} and the definition of $V^*$ because 
for each $\pi\in \Pi$, an unconstrained MDP problem $(X,A^{\pi}, P_{\pi},R_{\pi})$
(restricted to $A^{\pi}$)
can be formulated with $P_{\pi,xy}^a = P_{xy}^a$ and $R_{\pi}(x,a) = R(x,a)$ for
$x$ and $y$ in $X$ and $a\in A^{\pi}(x)$. The policy $\pi^*$ in the statement of the theorem is an optimal 
(or uniformly-optimal) policy
for the MDP and is in $F(\pi)$.
\endproof
\vspace{0.3cm}
Define $T_F:B(X\times \Pi) \rightarrow B(X\times \Pi)$ such that
for any $u\in B(X\times \Pi)$,
\[
   T_F(u)(x,\pi) = \max_{g\in F(\pi)} \left ( R(x,g(x)) + \gamma \sum_{y\in X} P_{xy}^{g(x)} u(y,g)\right ), x\in X, \pi\in \Pi.
\]
\begin{thm}
$T_F$ is a $\gamma$-contraction mapping on $B(X\times \Pi)$.
\end{thm}
\proof
For any $u$ and $v$ in $B(X\times \Pi)$,
\begin{eqnarray*}
\lefteqn{||T_F(u) - T_F(v)||_{\infty}}\\
& & \leq \max_{x\in X} \max_{\pi\in\Pi} \max_{g\in F(\pi)} \left ( \gamma \sum_{y\in X} P_{xy}^{g(x)} \left |u(y,g) - v(y,g) \right | \right ) \mbox{ by Prop A.3 in~\cite{hern}} \\
& & \leq \gamma ||u-v||_{\infty}
\end{eqnarray*}
\endproof

It follows that from the Banach fixed-point theorem, there exists $u\in B(X\times \Pi)$ that uniquely 
satisfies $T_F(u)=u$. We show below $u=V^*$.
The resulting optimality equation is essentially based on Chen and Blankenship's 
DP-equations restricted to $X\times \Pi$ with some proper modifications.
However, most importantly it works under Assumption~\ref{ass:ext}
and over $\Pi$ (not over the history-dependent policy set).
Furthermore, our proof does not use the limit arguments (with building up a 
history-dependent policy) for a sequence of finite-horizon values as in~\cite{chen}.

\vspace{0.2cm}
\begin{thm}
\label{thm:opteq}
$T_F(V^*) = V^*.$
\end{thm}
\vspace{0.2cm}

It \emph{is} possible that some $\pi$-infeasible policy $g$ in $\Pi\setminus \Phi(\pi)$ attains $V^*$, i.e., $V^* = V^g$, in which case $V^g$ also satisfies $T_F(V^g)=V^g$ from the uniqueness. Putting in another way, the value functions of some $\pi$-feasible policy and $\pi$-infeasible policy, respectively, can satisfy simultaneously the optimality equation.

\proof
Fix $x$ in $X$ and $\pi$ in $\Pi$.
We first show that $V^*(x,\pi) \leq T_F(V^*)(x,\pi)$.
\begin{eqnarray*}
\lefteqn{V^*(x,\pi) = \max_{g\in F(\pi)}V^g(x)  \mbox{ by the definition of } V^*}\\
& & = \max_{g\in F(\pi)} \left ( R(x,g(x)) + \gamma \sum_{y\in X} P_{xy}^{g(x)} V^g(y) \right ) \mbox{ by the definition of } V^g\\
& & \leq \max_{g\in F(\pi)} \left ( R(x,g(x)) + \gamma \sum_{y\in X} P_{xy}^{g(x)} 
\max_{\sigma\in F(g)} V^\sigma(y) \right ) \mbox{ because } g\in F(g) \\
& & = \max_{g\in F(\pi)} \left ( R(x,g(x)) + \gamma \sum_{y\in X} P_{xy}^{g(x)} V^*(y,g) \right ) \mbox{ by the definition of } V^*\\
& & = T_F(V^*)(x,\pi).
\end{eqnarray*} On the other hand, due to the existence of $\pi^*\in F(\pi)$ 
such that $V^*(y,\pi) = V^{\pi^*}(y)$ for all $y\in X$ by Theorem~\ref{thm:ext}, we have that 
\begin{eqnarray*}
\lefteqn{V^*(x,\pi) = V^{\pi^*}(x)}\\
& &   = R(x,\pi^*(x)) + \gamma \sum_{y\in X} P_{xy}^{\pi^*(x)} V^{\pi^*}(y)\\
& &   = \max_{g\in F(\pi)} \left( R(x,g(x)) + \gamma \sum_{y\in X} P_{xy}^{g(x)} V^{\pi^*}(y) \right ) \mbox{ by } T_F(V^*)=V^* \mbox{ and the property of } \pi^* \\
& &   \geq R(x,\sigma(x)) + \gamma \sum_{y\in X} P_{xy}^{\sigma(x)} V^{\pi^*}(y) \mbox{ for any } \sigma \in F(\pi) \\
& &   = R(x,\sigma(x)) + \gamma \sum_{y\in X} P_{xy}^{\sigma(x)} V^*(y,\pi) \mbox{ by the property of } \pi^* \\
& &   \geq R(x,\sigma(x)) + \gamma \sum_{y\in X} P_{xy}^{\sigma(x)} V^*(y,\sigma) \mbox{ because } \sigma \in F(\pi) 
\end{eqnarray*} Because the last inequality holds for any $\sigma$, it follows that $V^*(x,\pi) \geq T_F(V^*)(x,\pi)$.
This proves the theorem. 
\endproof

\begin{cor}
Given $\pi\in \Pi$,
consider $\phi \in \Pi$ such that for all $x\in X$, $\phi(x) \in A^*(x)\cap A^{\pi}(x)$, where
\[
 A^*(x) = \left\{ \sigma(x) \biggl | \sigma \in \argmax_{g\in F(\pi)} \left ( R(x,g(x)) + \gamma \sum_{y\in X} P_{xy}^{g(x)} V^*(y,g) \right ) \right \}, x\in X.
\] 
Then $\phi\in F(\pi)$ and $V^\phi(x)=V^*(x,\pi)$ for all $x\in X$.
\end{cor}

Because there exist $\pi^*\in F(\pi)$ such that $V^*(x,\pi) = V^{\pi^*}(x)$ for all $x\in X$, 
$A^*(x)$, $A^{\pi}(x)$, and their intersection are all non-empty for all $x\in X$.

\proof
Fix $x$ in $X$. By Theorem~\ref{thm:ext} and~\ref{thm:opteq}, for some $\pi^*$
\begin{eqnarray*}
\lefteqn{V^{\pi^*}(x) = \max_{g\in F(\pi)} \left( R(x,g(x)) + \gamma \sum_{y\in X} P_{xy}^{g(x)} V^{\pi^*}(y) \right )} \\
& &  = R(x,\phi(x)) + \gamma \sum_{y\in X} P_{xy}^{\phi(x)} V^{\pi^*}(y) \mbox{ by the definition of } \phi\\
& & = T_{\phi}(V^{\pi^*})(x).
\end{eqnarray*} Therefore, $V^{\pi^*} = T_{\phi}(V^{\pi^*})$, which implies that $V^{\phi}=V^{\pi^*}$ (see, e.g.,~\cite{ross}).
The feasibility follows from $\phi(x) \in A^{\pi}(x)$ for all $x\in X$. 
\endproof

\section{Computational Algorithms}
\label{sec:algo}
\subsection{Synchronous Off-line Method}

Recall that the goal is to (approximately) solve $M^c_{\pi^c}(x_0)$.
We start this subsection with providing a trivial PI-type algorithm
for finding a policy $\pi^*$ in $F(\pi^c)$ that achieves $V^{\pi^*}(x) = V^*(x,\pi^c)$ for all $x$ in $X$. This leads to a lower-bound solution such that $V^*(x_0,\pi^c) \leq V^*_c(x_0,\pi^c)$.
The algorithm is the same as PI but searches over $F(\pi^c)$ (instead of $\Pi$).
At each iteration $t\geq 1$,
\begin{eqnarray*}
\lefteqn{\pi_{t+1} \in \biggl \{ f \biggl | f\in F(\pi^c), \forall x\in X,
R(x,f(x)) + \gamma \sum_{y\in X} P_{xy}^{f(x)} V^{\pi_t}(y,\pi^c)}\\
& & \hspace{3cm} =
\max_{g\in F(\pi^c)} \biggl ( R(x,g(x)) + \gamma \sum_{y\in X} P_{xy}^{g(x)} V^{\pi_t}(y,\pi^c) \biggr ) \biggr  \},
\end{eqnarray*} where we write $V^{g}(x)$ as $V^{g}(x,\pi^c)$ (with the abuse of the notation) for $g \in F(\pi^c)$. 
The algorithm stops if $V^{\pi_{t+1}}(x,\pi^c) = V^{\pi_{t}}(x,\pi^c)$ for all $x\in X$ and starts with arbitrary policy in $F(\pi^c)$ as $\pi_1$. Because $\pi^c\in F(\pi^c)$, $\pi^c$ would be an immediate choice for $\pi_1$.

A VI-type algorithm is also trivially given from the property of $T_F$.
An interesting exact algorithm would be a \emph{backward induction}-type method (like solving a finite-horizon MDP) that works by recursively calling subproblems 
(at the root-$(x_0,\pi^c)$ level) based on the DP-equation until we hit the problems whose induced policy sets are singleton (at the leaf level).
In a backward or bottom-up fashion, once recursively called subproblems are solved, the uniformly-optimal value functions are returned to the problem that called the subproblems.

The question is whether we can improve $\pi^*$ obtained by such DP-based algorithms or not.
Unless $\Phi(\pi^c) = F(\pi^c)$, we can try searching further $\Phi(\pi^c)\setminus F(\pi^c)$ for a better policy 
than $\pi^*$ for $M^c_{\pi^c}(x_0)$. 
Unfortunately, the main obstacle would be about characterizing the search space that makes exploitable during a solution process. 
We circumvent this by the idea of generating a sequence of the sets induced from a \emph{unifomly-improving} sequence of the
feasible policies (as used before for related problems in~\cite{chang2}\cite{changexact}).
Unfortunately, (to the author's best knowledge) we still do not have an iterative method of policy-improvement ``focused at the state $x_0$," rather than uniformly,
for generating a sequence of $\pi^c$-feasible policies $\{\phi_t\}$ such that $V^{\phi_{t+1}}(x_0) \geq V^{\phi_{t}}(x_0)$ for $t\geq 1$ converging to $V^*(x_0,\pi^c)$
except for the just one iteration with updating at $x_0$.
Because still \emph{interdependencies} among the states for improvement need to be incorporated, developing such algorithm seems challenging.

For a given $\pi\in \Phi(\pi^c)$, let $\alpha_{\pi}$ be a subset of $A$ such that
\begin{equation}
\label{eqn:alpha}
\alpha_{\pi}(x) = \left \{a \biggl | a\in A(x), C(x,a) + \beta \sum_{y\in X} P_{xy}^a J^{\pi}(y) \leq J^{\pi}(x) + \Theta_{\pi}(x) \right \}, x\in X,
\end{equation} where $\Theta_{\pi}$ is selected from $B(X)$.
Also let $\Pi_{\pi}$ be $\{g | g\in \Pi, g(x)\in \alpha_{\pi}(x) \forall x\in X\}$ induced from $\alpha_{\pi}$.
It turns out that
for any $g\in \Pi_{\pi}$, $J^g(x) \leq J^{\pi}(x) + \Theta_{\pi}(x)/(1-\beta)$ for all $x\in X$~\cite{changimp}.

We consider the two extreme cases for simplicity in this paper.
If we just set $\Theta_{\pi}(x) = 0$ for all $x$ in $X$, $\alpha_{\pi} = A^{\pi}$, where $A^{\pi}$ is the element of the example DP-inducible set in~(\ref{eqn:DPset}). On the other hand, we can set
\[
\Theta_{\pi}(x) = (1-\beta)\left (J^{\pi^c}(x) - J^{\pi}(x) \right ), x\in X,
\] expressing the consumable \emph{slackness} for the feasibility of $\pi$ relative to $\pi^c$.
While the former case yields that $\Pi_{\pi} \subseteq \Phi(\pi)$,
the latter case does that $\Pi_{\pi} \subseteq \Phi(\pi^c)$.
Note that $\Pi_{f}$ does not necessarily include $\Pi(g)$ even if $g\in \Pi(f)$ even if
$\Phi(g) \subseteq \Phi(f)$ if $g \in \Phi(f)$.
In other words, the subset relationship does not hold necessarily for the policy sets obtained by the first case and the second case, respectively.

We have 
the following iterative algorithm $\mathcal{A}$ for possibly improving $\pi^*$: Set $\pi_1=\pi^*$. For $t\geq 1$,
\[
    \pi_{t+1} \in \left \{\pi \biggl | \pi\in \Pi_{\pi_t}, V^{\pi}(x) = \max_{\phi\in \Pi_{\pi_t}} V^{\phi}(x) \mbox{ for all } x\in X\right \},
\] where we stop the iteration if $V^{\pi_{t+1}}=V^{\pi_t}$ and $J^{\pi_{t+1}}=J^{\pi_t}$
and $\alpha_{\pi_{t+1}}=\alpha_{\pi_t}$.

First, obtaining the set for $\pi_{t+1}$ corresponds to
solving an MDP $(X,\alpha_{\pi_t},P_{\pi_t},R_{\pi_t})$ where $P_{\pi_t}$ and $R_{\pi_t}$
are obtained from $P$ and $R$, respectively, by restricting to $\alpha_{\pi_t}$.
That is, $\pi_{t+1}$ is uniformly-optimal among the policies in $\Pi_{\pi_t}$.

Second, by an inductive argument on $t$, for any $t\geq 1$ and $x$ in $X$,
\[
   V^{\pi^*}(x) \leq \max_{g \in \bigcup_{k=1}^{t} \Pi_{\pi_k}} V^g(x) \leq V^{\pi_{t+1}}(x) \leq V^*_c(x,\pi^c)
\] In other words, the sequence of the value functions of $\{\pi_t\}$ is monotonically and uniformly improving toward $V^*_c$ leading to a local uniform-optimality.
If $\pi^c_*$ is an element of any $\Pi_{\pi_t}$, we have a global convergence.

Lastly, more careful attention needs to be paid to the stopping condition than
to the methods used in~\cite{chang2}~\cite{changexact}.
Because $\Pi$ is finite and the monotonicity of $V^{\pi_{t+1}} \geq V^{\pi_t}$ holds for all $t$, 
at some $t$, we must have $V^{\pi_{t+1}} = V^{\pi_t}$. In this case, both $\pi_{t+1}$ 
and $\pi_t$ are uniformly-optimal policies for the MDP $(X,\alpha_{\pi_t},P_{\pi_t},R_{\pi_t})$.
This is obvious
because $V^{\pi_t}(x) = R(x,\pi_{t+1}(x)) + \gamma \sum_{y\in X} P_{xy}^{\pi_{t+1}(x)} V^{\pi_{t+1}}(y)
=\max_{a\in \alpha_{\pi_t}(x)}( R(x,a)+ \gamma \sum_{y\in X} P_{xy}^a V^{\pi_t}(y))$ for all $x$ in $X$. 
Now, if $\pi_t=\pi_{t+1}$, then we stop. But if not, it is \emph{still} possible to improve upon $\pi_{t+1}$.
We need to check further whether $J^{\pi_{t+1}} = J^{\pi_t}$ or not.
If there exists $x$ in $X$ such that $J^{\pi_{t+1}}(x) \neq J^{\pi_t}(x)$, then $\alpha_{\pi_{t+2}}$ can be different from $\alpha_{\pi_{t+1}}$.
The general stop condition is thus when $V^{\pi_{t+1}}=V^{\pi_t}$ and $J^{\pi_{t+1}}=J^{\pi_t}$
and $\alpha_{\pi_{t+1}}=\alpha_{\pi_t}$ so that no more changes can be made after this happens.

It is remarkable that at each iteration of $\mathcal{A}$, an MDP associated with $\alpha_{\pi_t}$ is solved and this is the key which leads to the monotonicity while keeping each policy $\pi_t$ generated \emph{feasible.} Designing such an algorithm is not straightforward.
The main difficulty lies on how to combine multiple feasible-policies to generate a feasible policy while uniformly improving the reward value-functions of the given policies.
Suppose that
a uniformly-optimal $\pi^c$-feasible policy $\sigma^*$ in $\Pi_{\sigma}$ and a uniformly-optimal $\pi^c$-feasible policy $\rho^*$ in $\Pi_{\rho}$ have been obtained by
solving $(X,\alpha_{\sigma},P_{\sigma},R_{\sigma})$ and $(X,\alpha_{\rho},P_{\rho},R_{\rho})$ for some $\sigma$ and $\rho$ in $\Phi(\pi^c)$, respectively. 
The two policies $\sigma^*$ and $\rho^*$ are \emph{incomparable} in general.
In other words, it is possible that neither $V^{\sigma^*} \geq V^{\rho^*}$ nor 
$V^{\sigma^*} \geq V^{\rho^*}$. 
In order to have a monotonicity, we need a policy that uniformly improves both policies at all states. A possible way is to use policy switching or parallel rollout~\cite{changbook} to ``mix" $\sigma^*$ and $\rho^*$ into a new policy $\psi$. 
Even if it results that $V^{\psi} \geq V^{\rho^*}$ and $V^{\psi} \geq V^{\sigma^*}$,
$\psi$ is \emph{not necessarily} $\pi^c$-feasible.
Indeed, if the open problem of multi-policy improvement were resolved, the following method is naturally drawn for approximately solving $M^c_{\pi^c}(x_0)$.
We generate a sequence of $\{\pi_k\}$ by applying PI to 
the MDP $(X,A,P,C)$ with $\pi_1=\pi^c$. We have that $J^{\pi^c} \geq J^{\pi_2} \geq J^{\pi_3} \cdots \geq J^{\pi_n}$ for some $n$.
Because $J^{\pi^*_c}$ must reside between the value functions of some consecutive policies, 
it is likely (at least promising) for $\bigcup_{t=1}^n \Pi_{\pi_t}$ to cover or approximate 
well $\Phi(\pi^c)$. A feasible policy $\pi^c_{\app}$ that improves all policies 
in $\bigcup_{t=1}^n \Pi_{\pi_t}$ can be obtained.
If $\pi^*$ is available, the same process is done with $\pi_1=\pi^*$
obtaining $\pi^*_{\app}$. $\pi^c_{\app}$ and $\pi^*_{\app}$ can be mixed as an approximate solution to $M^c_{\pi^c}(x_0)$.

The next question is if we can still make further improvement once $\mathcal{A}$ converges.
We end this subsection with a brief discussion about a possible way.
Suppose that the algorithm $\mathcal{A}$ stopped at $t=n$. 
Now consider applying a single policy improvement of PI to $\pi_n$ within the MDP
$(X,A,P,R)$ \emph{with the original parameters}. That is,
let $\PI(\pi_n)\in \Pi$ be given as 
\[
\PI(\pi_n)(x) \in \argmax_{a\in A(x)} \left ( R(x,a) + \gamma \sum_{y\in X}P_{xy}^a V^{\pi_n}(y) \right ), x\in X.
\] 
If $\PI(\pi_n)$ is in $\Phi(\pi^c)$ and
$V^{\PI(\pi_n)} = V^{\pi_n}$, it follows that $V^{*}_{(X,A,P,R)}(x_0) = V^{\PI(\pi_n)}(x_0) \geq V^{\pi^*_c}(x_0) \geq V^{\pi_n}(x_0) = V^{\PI(\pi_n)}(x_0)$, where $V^*_{(X,A,P,R)}(x):= \max_{\pi\in \Pi} V^{\pi}(x), x\in X,$ and $\pi^*_c$ is again a solution to $M^c_{\pi^c}(x_0)$.
We have a global convergence: $\pi_n$ is a solution to $M^c_{\pi^c}(x_0)$.
On the other hand,
if $\PI(\pi_n)$ is in $\Phi(\pi^c)$ but there exists $x$ in $X$ such that 
$V^{\PI(\pi_n)}(x) > V^{\pi_n}(x)$, then $\PI(\pi_n)$
is a \emph{strictly-improving} policy relative to $\pi_n$ for $M^c_{\pi^c}(x_0)$.
When we have this case or the case where $\PI(\pi_n)$ is not in $\Phi(\pi^c)$, we go to the next round of checking a possible improvement of $\pi_n$ with $\PI^2(\pi_n) = \PI(\PI(\pi_n))$.
This same process can continue until (the value functions of) the two consecutive policies are equal (as in the PI stop condition).

\subsection{Asynchronous On-line Method}

In this subsection, we interpret $t\geq 0$ as \emph{the system time}. The algorithm
below applies to the system in \emph{on-line} updating asynchronously \emph{only at the current state}.

At $t=0$, set $\pi_0 = \pi^c$ or $\pi^*$ obtained by off-line DP or any known feasible policy in $\Phi(\pi^c)$. The initial state is $x_0$.
At $t\geq 0$, we first update $\alpha_{\pi_{t}}$ \emph{only at} $x_t$ such that $\alpha_{\pi_{t}}(x)$ does not change for all $x\neq x_{t}$ and
\[
   \alpha_{\pi_t}(x_t) = \left \{ a \biggl | a\in A(x_t), C(x_t,a) + \beta \sum_{y\in X} P_{x_ty}^a J^{\pi_t}(y) \leq J^{\pi_t}(x_t) \right \}.
\]
Then $\pi_t$ is also updated \emph{only at} $x_t$: for $x\neq x_t$, $\pi_{t+1}(x) = \pi_t(x)$ and 
\[
   \pi_{t+1}(x_t) \in \argmax_{a\in \alpha_{\pi_t}(x_t)}\left ( R(x_t,a) + \gamma \sum_{y\in X} P_{x_ty}^a V^{\pi_t}(y) \right ).
\]
Once $\pi_{t+1}$ is obtained, $\pi_{t+1}(x_t)$ is taken and $x_t$ makes a random transition to $x_{t+1}$.

We claim that \emph{for all $t\geq 0$, $\pi_{t+1}$ is feasible}. 
By our choice, $\pi_0$ is feasible.
At $x\neq x_t$,
$U_{\pi_{t+1}}(J^{\pi_{t}})(x) 
= C(x,\pi_{t+1}(x)) + \beta \sum_{y\in X} P_{xy}^{\pi_{t+1}(x)} J^{\pi_t}(y) 
= C(x,\pi_{t}(x)) + \beta \sum_{y\in X} P_{xy}^{\pi_t(x)} J^{\pi_t}(y) 
= J^{\pi_t}(x)$.
On the other hand, at $x = x_t$, because $\pi_{t+1}(x_t)$ is in $\alpha_{\pi_{t}}(x_t)$,
$U_{\pi_{t+1}}(J^{\pi_{t}})(x) 
= C(x,\pi_{t+1}(x)) + \beta \sum_{y\in X} P_{xy}^{\pi_{t+1}(x)} J^{\pi_t}(y) 
\leq J^{\pi_t}(x)$. 
Putting together, $U_{\pi_{t+1}}(J^{\pi_t}) \leq J^{\pi_t}$ so that $J^{\pi_{t+1}}\leq J^{\pi_t}$.
Because $\pi_t$ is feasible, it follows that $\pi_{t+1}$ is feasible. The inductive argument on $t$ proves
the claim.
Furthermore, \emph{the sequence of the value functions of $\{\pi_t\}$ is monotonically and uniformly improving.}
At $x\neq x_t$, 
$T_{\pi_{t+1}}(V^{\pi_{t}})(x) 
= R(x,\pi_{t+1}(x)) + \gamma \sum_{y\in X} P_{xy}^{\pi_{t+1}(x)} V^{\pi_t}(y) 
= R(x,\pi_{t}(x)) + \gamma \sum_{y\in X} P_{xy}^{\pi_{t}(x)} V^{\pi_t}(y) 
= V^{\pi_t}(x)$.
At $x = x_t$, from the definition of $\pi_{t+1}(x_t)$,
$T_{\pi_{t+1}}(V^{\pi_{t}})(x) 
= R(x,\pi_{t+1}(x)) + \gamma \sum_{y\in X} P_{xy}^{\pi_{t+1}(x)} V^{\pi_t}(y) 
\geq V^{\pi_t}(x)$. It follows that $V^{\pi_{t+1}} \geq V^{\pi_t}$.

The convergence behaviour of the algorithm output 
here is closely related to that of the on-line algorithm studied in~\cite{changonline} 
for MDPs (cf., Section IV in~\cite{changonline}).
We omit a detailed discussion but end this section with an intuitive remark.
Assume that the underlying MDP for $M^c_{\pi^c}(x_0)$ is communicating. In other words, every Markov chain 
induced by fixing each policy in $\Pi$ is communicating. This assumption implies that
every state $x$ in $X$ is visited infinitely often within $\{x_t\}$. 
Because $\Pi$ is finite, it follows that
there exists a finite-time $K < \infty$ such that
for all $t > K$ $\pi_t$ is a uniformly-optimal $\pi^c$-feasible policy for the MDP $(X,\alpha_{\pi_K},P_{\pi_K},R_{\pi_K})$
and $V^{\pi_{t+1}}(x) \geq \max_{g \in \bigcup_{k=1}^{K} \Pi_{\pi_k}} V^g(x)$ for all
$x$ in $X$ achieving a local uniform-optimality for the CMDP problem $M^c_{\pi^c}(x_0)$.

\section{Concluding Remarks}

In the problem $M^c_{\pi^c}(x_0)$, the feasibility of a policy is determined by the value function of $\pi^c$. 
If the problem is set up with a threshold function $\kappa$ in $B(X)$ for the feasibility,
one can approximate $\kappa$, for example, by the value function of a policy in $\argmin_{\pi\in\Pi_{\kappa}} \sum_{x\in X} |\kappa(x) - J^{\pi}(x)|$, where $\Pi_{\kappa} = \{ g | g\in \Pi, g(x) \in A^{\kappa}(x) \mbox{ } \forall x\in X\}$, where
\[
 A^{\kappa}(x) = \left \{a \biggl | a\in A(x), C(x,a) + \beta \sum_{y\in X} P_{xy}^a \kappa(y) > \kappa(x) \right \}, x\in X.
\] An MDP $(X,A^{\kappa},P_{\kappa},C_{\kappa})$ restricted to $A^{\kappa}$ for $P$ and $C$ can be solved then to obtain $\pi$ that minimizes the distance with $\kappa$.

While discussing $\alpha_{\pi}$ in~(\ref{eqn:alpha}), the two possibilities were considered. The second case is not suitable for an example of $A^{\pi}$ in the DP-inducible set.
The issue is that as seen, for all $\pi\in \Pi$, we have that $F(\pi) \subseteq \Phi(\pi^c)$ in comparison with $F(\pi)\subseteq \Phi(\pi)$. 
We would lose the fundamental DP spirit because 
the substructure relationship of DP that a problem recursively calls its ``smaller" subproblems of having the same form with the smaller parameters inside the problem does not hold necessarily.

Even if another DP-inducible set could be found, finding $A^{\pi}$ such that the $\pi$-feasible policy set induced from $A^{\pi}$ provides a ``maximal" DP-inducible set such that $\Phi(\pi) \setminus F(\pi)$ is as small as possible for all $\pi$ in $\Pi$ seems very difficult.
Even finding an DP-inducible set $\{B^{\pi}\}$ that ``expands" $\{A^{\pi}\}$ in~(\ref{eqn:DPset}) such that $B^{\pi} \supset A^{\pi}$ for all $\pi\in \Pi$ is an open problem.

The DP-equation in Theorem~\ref{thm:opteq} is a generalization of Bellman's optimality equation because a finite-MDP problem is a special case of the CMDP problem solved by the DP-equation in Section~\ref{sec:dp}.

Even though the values (i.e., performances) of randomized policies can be (much) bigger than those of pure policies~\cite{feinberg}~\cite{frid} in infinite-horizon discounted CMDP problems,
besides representational and computational complexity-issues while realizing a randomized policy in actual implementation~\cite{chen}, the question of whether such CMDP problems as $M^c_{\pi^c}(x_0)$ have (at least approximate) DP characteristics or not must be an important issue to be addressed,
especially, in the pure, Markovian, and stationary policy set.
This is because any CMDP is basically an MDP and the fundamental theory of infinite-horizon discounted MDP is built on the DP-concept from the optimality principle of Bellman.
The DP-equation for MDP expresses the principle and the existence of an optimal policy in the policy set that attains the optimal values at all states follows from the equation.
This work shows that such a CMDP problem in general has an implicit DP-structure among the CMDP problems inducible from the original problem in the policy space, which can be exploited during a solution process.

\end{document}